\documentclass[11pt]{article}
\usepackage[fleqn]{amsmath}
\usepackage{amsfonts}

\newtheorem{theorem}{Theorem}[section]

\newtheorem{proposition}[theorem]{Proposition}

\DeclareMathOperator{\xGL}{GL}
\DeclareMathOperator{\xNS}{NS}
\DeclareMathOperator{\tr}{tr}
\DeclareMathOperator{\Kahler}{K\ddot{a}hler}

\newcounter{cnt1}
\newcounter{cnt2}

\begin{document}
\title{Salem Numbers and Automorphisms of Complex Surfaces}
\author{Paul Reschke}
\date{}

\maketitle

\begin{abstract}
For two-dimensional complex tori, we characterize the set of all values of positive entropy that arise from automorphisms. For K3 surfaces, we give sufficient conditions for a positive value to be the entropy of some automorphism.
\end{abstract}

\section{Overview}

In this note, we examine complex surface automorphisms with positive entropy. We focus on automorphisms of complex tori and K3 surfaces.
\\
\\
For complex tori, we determine precisely which positive values arise as entropies of such automorphisms. Each of these values is the natural logarithm of a Salem number of degree two, four, or six. (See \S 2 in this note.) We specify necessary and sufficient conditions for a Salem number of degree four or six to give the entropy of some torus automorphism, and we show that every Salem number of degree two gives the entropy of some torus automorphism.
\begin{theorem}
Let \(S(t)\) be the minimal polynomial for a Salem number \(\lambda\).
\begin{list}{\emph{\textbf{\arabic{cnt1})}}}{\setlength\leftmargin{24pt} \usecounter{cnt1}}
\item If \(\lambda\) has degree six, then \(\log(\lambda)\) is the entropy of some automorphism of a two-dimensional complex torus if and only if \(S(1)=-m^2\) for some integer \(m\) and \(S(-1)=n^2\) for some integer \(n\).
\item If \(\lambda\) has degree four, then \(\log(\lambda)\) is the entropy of some automorphism of a two-dimensional complex torus if and only if one of the following three cases holds: (a) \(S(1)=-m^2\) for some integer \(m\); (b) \(S(-1)=n^2\) for some integer \(n\); or (c) \(S(1)=-\frac{1}{2}m^2\) and \(S(-1)=\frac{1}{2}n^2\) for some integers \(m\) and \(n\).
\item If \(\lambda\) has degree two, then \(\log(\lambda)\) is the entropy of some automorphism of a two-dimensional complex torus.
\end{list}
These cases constitute all possible positive entropies for an automorphism of a two-dimensional complex torus.
\end{theorem}
(See \S 3 in this note.) This theorem extends results in \cite{gmc}, \cite{mc1}, and \cite{mc2} which provide specific examples of Salem numbers whose logarithms are entropies of automorphisms of two-dimensional complex tori. We discuss some additional results specific to automorphisms of projective tori and tori which are products of one-dimensional complex tori. (See \S 5 in this note.)
\\
\\
We consider an analogous investigation of automorphisms of K3 surfaces. The entropy of any such automorphism must be zero or the natural logarithm of a Salem number of degree at most twenty-two. Kummer surfaces (K3 surfaces arising from quotients of tori) have automorphisms whose entropies are natural logarithms of Salem numbers of degree at most six. Indeed, any Salem number that gives the entropy of a torus automorphism must also give the entropy of the corresponding Kummer surface automorphism. We identify a sufficient condition for a Salem number of degree fourteen to give the entropy of some automorphism of a K3 surface. A result in \cite{gmc} is that the same condition is sufficient for a Salem number of degree twenty-two to give the entropy of some K3 surface automorphism.
\begin{theorem}
Let \(S(t)\) be the minimal polynomial for a Salem number \(\lambda\) of degree fourteen satisfying \(S(-1)S(1)=-1\). Then \(\log(\lambda)\) is the entropy of some K3 surface automorphism. 
\end{theorem}
(See \S 6 in this note.) This theorem comes from combining results in \cite{gmc}, \cite{mc2}, and \cite{ogu}. We note that there are examples of K3 surface automorphisms with positive entropy not accounted for by either this theorem or the Kummer construction. (See \cite{mc2}, and \S 6 in this note.)
\\
\\
If a compact \(\Kahler\) surface admits an automorphism with positive entropy, then (after any exceptional curves with finite order under the automorphism are contracted) the surface can only be a complex torus, a K3 surface, an Enriques surface, or a rational surface. (See \cite{can}, \S 2.) The set of possible entropies for rational surface automorphisms is described in \cite{ueh}. 

\section{Computing Entropies on a Torus}

Let \(X=\mathbb{C}^2/\Lambda\) be a two-dimensional complex torus; so \(\Lambda \subseteq \mathbb{C}^2\) is a four-dimensional lattice acting holomorphically on \(\mathbb{C}^2\) by translations. If \(F\) is an automorphism of \(X\), then \(F\) is necessarily the quotient of some invertible linear map \(\widetilde{F} \in \xGL_2(\mathbb{C})\) with \(\widetilde{F}(\Lambda)=\Lambda\). In this sense, \(F\) can be expressed as a \(2 \times 2\) matrix over \(\mathbb{C}\). 
\\
\\
Given coordinates \(z_1\) and \(z_2\) for each factor of \(\mathbb{C}^2\), the one-forms \(dz_1\), \(dz_2\), \(d\overline{z}_1\), and \(d\overline{z}_2\) on \(\mathbb{C}^2\) are all invariant under translations and hence descend to one-forms on \(X\). Moreover, these four one-forms generate the cohomology ring for \(X\) with coefficients in \(\mathbb{C}\). If \(F\) is an automorphism of \(X\), then \(F\) induces an automorphism \(F^*\) of \(H^*(X;\mathbb{C})\) via the pull-back action. As an element of \(\xGL_4(\mathbb{C})\), the pull-back action of \(F\) on \(H^1(X;\mathbb{C})\) is the direct sum of the transpose of \(F\) and the conjugate transpose of \(F\); that is, \(F^*=F^T \oplus \bar{F}^T\). The pull-back action of \(F\) gives a convenient means of computing the entropy of \(F\) via the spectral radius \(\rho(F^*)\). Indeed, the entropy of \(F\) is given by \[h(F)=\log(\rho(F^*:H^2(X;\mathbb{C}) \rightarrow H^2(X;\mathbb{C}))).\]
This formula is a special case of a general statement for automorphisms of compact \(\Kahler\) manifolds due to Gromov and Yomdin. (See \cite{can}, \S 2.1.)
\\
\\
If the eigenvalues of \(F\) are \(\gamma_1\) and \(\gamma_2\), then the eigenvalues of \(F^*\) restricted to \(H^1(X;\mathbb{C})\) are \(\gamma_1\), \(\gamma_2\), \(\overline{\gamma_1}\), and \(\overline{\gamma_2}\). Since the pull-back commutes with the wedge product, the eigenvalues of \(F^*\) restricted to \(H^2(X;\mathbb{C})\) are \(|\gamma_1|^2\), \(|\gamma_2|^2\), \(\gamma_1 \gamma_2\), \(\gamma_1 \overline{\gamma_2}\), \(\overline{\gamma_1} \gamma_2\), and \(\overline{\gamma_1} \overline{\gamma_2}\). Thus the entropy of \(F\) is given by
\[h(F)=\log(\max\{|\gamma_1|^2,|\gamma_2|^2\}).\]
\\
Since \(F^*\) preserves the lattice \(H^1(X;\mathbb{Z})\) inside \(H^1(X;\mathbb{C})\), \(F^*\) must be an element of \(\xGL_4(\mathbb{Z})\) and hence must have determinant one or negative one; so \(|\gamma_2|=|\gamma_1|^{-1}\) and the eigenvalues of \(F^*\) restricted to \(H^2(X;\mathbb{C})\) are some \(\lambda \geq 1\), \(\lambda^{-1}\), \(\alpha\), \(\overline{\alpha}\), \(\beta\), and \(\overline{\beta}\) with \(|\alpha|=|\beta|=1\). The entropy of \(F\) is then \(\log(\lambda)\). Since \(F^*\) also preserves the lattice \(H^2(X;\mathbb{Z})\) inside \(H^2(X;\mathbb{C})\), the characteristic polynomial for \(F^*\) restricted to \(H^2(X;\mathbb{C})\) must be some monic reciprocal polynomial \(Q(t) \in \mathbb{Z}[t]\) of degree six. The minimal polynomial for \(\lambda\) is then some irreducible factor of \(Q(t)\); if \(\lambda\) is greater than one, then this irreducible factor is a reciprocal polynomial with exactly two roots off the unit circle, both of which are positive. Thus \(\lambda\) is either one or a Salem number of degree two, four, or six. (See also \cite{mc1}, \S 3.)

\section{Constructing Torus Automorphisms}

Let \(F\) be an automorphism of a two-dimensional complex torus \(X\) with eigenvalues \(\gamma_1\) and \(\gamma_2\). Let \(Q(t) \in \mathbb{Z}[t]\) and \(P(t) \in \mathbb{Z}[t]\) be the characteristic polynomials of \(F^*\) restricted to \(H^2(X;\mathbb{C})\) and \(H^1(X;\mathbb{C})\), respectively. Then
\[Q(t)=(t-|\gamma_1|^2)(t-|\gamma_2|^2)(t-\gamma_1\gamma_2)(t-\gamma_1\overline{\gamma_2})(t-\overline{\gamma_1}\gamma_2)(t-\overline{\gamma_1}\overline{\gamma_2})\]
\[\indent =t^6+at^5+bt^4+ct^3+bt^2+at+1\]
and
\[P(t)=(t-\gamma_1)(t-\gamma_2)(t-\overline{\gamma_1})(t-\overline{\gamma_2})\]
\[\indent =t^4+jt^3-at^2+kt+1\]
for some integers \(a\), \(b\), \(c\), \(j\), and \(k\) with \(jk=b+1\) and \(j^2+k^2=-c-2a\). So \(Q(1)=-(j-k)^2=-m^2\) for some integer \(m\) and \(Q(-1)=(j+k)^2=n^2\) for some integer \(n\); we say that \(Q(t)\) has the ``square property''. (Compare \cite{gmc}, \S 1.)
\\
\\
In particular, any Salem number of degree six that gives the entropy for some torus automorphism must have a minimal polynomial with the square property. The following existence theorem provides a means of determining precisely when a Salem number gives the entropy of some torus automorphism.
\begin{theorem}
[\cite{mc1}, Theorem 4.4] Let \(P(t) \in \mathbb{Z}[t]\) be monic of degree four with roots \(\gamma_1\), \(\gamma_2\), \(\overline{\gamma_1}\), and \(\overline{\gamma_2}\) such that \(|\gamma_2|=|\gamma_1|^{-1}\). Then there is some two-dimensional complex torus \(X\) and some automorphism \(F\) of \(X\) such that \(P(t)\) is the characteristic polynomial of \(F^*\) restricted to \(H^1(X;\mathbb{C})\).
\end{theorem}
Now let \(Q(t) \in \mathbb{Z}[t]\) be a monic reciprocal polynomial of degree six with the square property; so
\[Q(t)=t^6+at^5+bt^4+ct^3+bt^2+at+1\]
for some integers \(a\), \(b\), and \(c\), and
\[Q(1)=2+2a+2b+c=-m^2\]
for some integer \(m\), and
\[Q(-1)=2-2a+2b-c=n^2\]
for some integer \(n\). Then \(m\) and \(n\) are either both odd or both even, and thus \(j=\frac{1}{2}(m+n)\) and \(k=\frac{1}{2}(n-m)\) are both integers. Since \(jk=b+1\) and \(j^2+k^2=-c-2a\), the roots of \(Q(t)\) must be the products of distinct roots of the polynomial
\[P(t)=t^4+jt^3-at^2+kt+1.\]
If \(Q(t)\) has exactly two roots off the unit circle, both positive, then the roots of \(P(t)\) occur in conjugate pairs. These calculations yield the following proposition.
\begin{proposition}
Let \(Q(t) \in \mathbb{Z}[t]\) be monic and reciprocal of degree six with exactly two roots off the unit circle, both positive. Then the following two statements are equivalent:
\begin{list}{\emph{\textbf{\arabic{cnt2})}}}{\setlength\leftmargin{24pt} \usecounter{cnt2}}
\item there is a two-dimensional complex torus \(X\) and an automorphism \(F\) of \(X\) such that \(F^*\) restricted to \(H^2(X;\mathbb{C})\) has characteristic polynomial \(Q(t)\); and
\item \(Q(t)\) has the square property.
\end{list}
\end{proposition}
Since the entropy of a torus automorphism must be zero or a Salem number, no polynomial with more than two roots off the unit circle can be the characteristic polynomial for the induced action on second cohomology. If some automorphism gives rise to a characteristic polynomial with no roots off the unit circle, then the automorphism must have entropy zero.
\\
\\
In particular, if \(\lambda\) is a Salem number of degree six, then \(\log(\lambda)\) is the entropy of some torus automorphism if and only if the minimal polynomial for \(\lambda\) has the square property. Similarly, a Salem number of degree two or four gives the entropy of some torus automorphism if and only if its minimal polynomial is a factor of some monic reciprocal polynomial of degree six with the square property and only two roots off the unit circle.
\\
\\
\emph{Proof of Theorem 1.1.} The degree-six case is evident from Proposition 3.2. In the degree-two case, if \(S(t)\) is the minimal polynomial for \(\lambda\), then \(S(t)(t-1)^2(t+1)^2\) is a monic reciprocal polynomial of degree six with only two roots off the unit circle satisfying \(S(1)=S(-1)=0\). In the degree-four case, if \(S(t)\) is the minimal polynomial for \(\lambda\), then \(S(t)C(t)\) is a monic reciprocal polynomial of degree six with only two roots off the unit circle if and only if \(C(t)\) is one of \(t^2-2t+1\), \(t^2-t+1\), \(t^2+1\), \(t^2+t+1\), or \(t^2+2t+1\); furthermore, \(S(t)C(t)\) has the square property for one of these choices of \(C(t)\) if and only if one of the three stated cases holds for \(S(t)\). \(\Box\)
\\
\\
If \(F\) is a torus automorphism with positive entropy, then \(F\) must have eigenvalues \(\gamma_1\) and \(\gamma_2\) with \(|\gamma_2|=|\gamma_1|^{-1} \neq 1\). So the eigenvalues of \(F^*\) restricted to \(H^1(X;\mathbb{C})\) cannot have magnitude one and in particular cannot be roots of unity. This result implies that \(F\) must be ergodic. (See \cite{mc1}, \S 4, and \cite{man}, III, Theorems 3.1, 8.5, and 8.6.)

\section{Lattice Isometries}

Given a two-dimensional complex torus \(X\), the intersection form makes \(H^2(X;\mathbb{Z})\) into an even unimodular lattice of signature (3,3), which we denote \(L_X\). If \(F\) is an automorphism of \(X\), then \(F^*\) is an isometry from \(L_X\) to itself. The isometry \(F^*\) extends to an isometry on \(L_X \otimes \mathbb{Q}\) and hence can be expressed in rational canonical form via some change of basis over \(\mathbb{Q}\). The invariant factors of \(F^*\) correspond to \(F^*\)-invariant subspaces of \(L_X \otimes \mathbb{Q}\); the same is true for any irreducible factors with multiplicity one within an invariant factor of \(F^*\). Since any subspace of \(L_X \otimes \mathbb{Q}\) must intersect \(L_X\) itself, any \(F^*\)-invariant subspace of \(L_X \otimes \mathbb{Q}\) must contain some \(F^*\)-invariant sublattice of \(L_X\). Moreover, since classes in \(H^2(X;\mathbb{Z})\) can be represented as isotopy classes of smoothly embedded real surfaces in \(X\) (whose intersection numbers are compatible with the intersection form on \(L_X\)), \(F^*\)-invariant sublattices of \(L_X\) must correspond to \(F\)-invariant spaces of isotopy classes. (See \cite{gom}, \S 1.2, and \cite{mc2}, \S 2.)
\\
\\
The study of torus automorphisms fits as a special case into the study of isometries of even unimodular lattices in general. Results in \cite{gmc} give a sufficient condition for a Salem number to be the entropy of an isometry of an even unimodular lattice: let \(\lambda\) be a Salem number with minimal polynomial \(S(t)\) of degree \(d\) such that
\[d \equiv 2 \mod{4},\]
and let \(p\) and \(q\) be positive integers satisfying \(p+q=d\) and
\[p \equiv q \mod{8};\]
then \(\log(\lambda)\) is the entropy of some isometry of the even unimodular lattice of signature (\(p\),\(q\)) if
\(S(-1)S(1)=-1\). (See \cite{gmc}, \S 1, and \cite{boy}, \S 2.) Moreover, a necessary condition for \(\log(\lambda)\) to be the entropy of such an isometry is that \(|S(-1)|\), \(|S(1)|\), and \(-S(-1)S(1)\) are all squares of integers. Let \(L_{3,3}\) be the (unique up to isometry) even unimodular lattice of signature (3,3). In this case, the classification of Salem numbers with respect to entropies of torus automorphisms gives a complete statement: let \(S(t)\) be the minimal polynomial for a Salem number of degree six; then \(S(t)\) is the characteristic polynomial for some isometry of \(L_{3,3}\) if and only if \(S(t)\) has the square property. (For a degree-six Salem polynomial \(S(t)\), the condition \(S(-1)S(1)=-1\) forces \(S(-1)=-1\) and \(S(1)=1\).)

\section{Projective Tori}

Any two-dimensional complex torus \(X\) has as part of its Hodge structure
\[H^2(X;\mathbb{C}) \cong H^{2,0} \oplus H^{1,1} \oplus H^{0,2},\]
and the intersection form extends to a Hermitian inner product on \(H^2(X;\mathbb{C})\) with signature (2,0) on \(H^{2,0} \oplus H^{0,2}\) and signature (1,3) on \(H^{1,1}\). If \(F\) is an automorphism of \(X\), then \(F^*\) must preserve the Hodge structure and the N\'eron-Severi group
\[\xNS(X) = H^{1,1} \cap H^2(X;\mathbb{Z}).\]
The torus \(X\) is projective if and only if \(\xNS(X)\) contains some element with positive self-intersection. So \(F\) can only be an automorphism of a projective torus if \(F^*\) preserves a sublattice of \(H^2(X;\mathbb{Z})\) of rank at least one and at most four. In particular, the characteristic polynomial of \(F^*\) cannot be irreducible if \(X\) is projective. Thus no Salem number of degree six can give the entropy of a projective torus automorphism. The projective torus \(\mathbb{C}/\mathbb{Z}[\zeta_3] \times \mathbb{C}/\mathbb{Z}[\zeta_3]\), where \(\zeta_3\) is a third root of unity, does admit an automorphism whose entropy is the logarithm of a Salem number of degree four (the smallest Salem number of degree four, in fact). (See \cite{mc2}, \S 1.)
\\
\\
Any complex torus of the form \(E \times E\), where \(E\) is a one-dimensional complex torus, is necessarily projective. Any linear map \(A \in \xGL_2(\mathbb{Z})\) with \(\det(A)=1\) gives an automorphism of \(E \times E\) with eigenvalues \(\frac{1}{2}(\tr(A)+\sqrt{\tr(A)^2-4})\) and \(\frac{1}{2}(\tr(A)-\sqrt{\tr(A)^2-4})\). (See \cite{mc1}, \S 4.) The characteristic polynomial for \(A^*\) restricted to \(H^2(E \times E;\mathbb{C})\) is then
\[(t^2-(\tr(A)^2-2)t+1)(t-1)^4,\]
and \(A\) has positive entropy if \(|\tr(A)|>2\). Similarly, if \(\det(A)=-1\), then the characteristic polynomial for \(A^*\) restricted to \(H^2(E \times E; \mathbb{C})\) is
\[(t^2-(\tr(A)^2+2)t+1)(t+1)^4,\]
and \(A\) has positive entropy if \(|\tr(A)|>0\). These calculations give sufficient conditions for a Salem number of degree two to give the entropy of an automorphism of \(E \times E\): let \(\lambda\) be a Salem number of degree two and let \(t^2-at+1\) (with \(a>2\)) be its minimal polynomial; then \(\log(\lambda)\) is the entropy of some automorphism of \(E \times E\) if \(a=m^2+2\) for some integer \(m\) or \(a=m^2-2\) for some integer \(m\). These examples do not include the case referred to above of an automorphism of \(\mathbb{C}/\mathbb{Z}[\zeta_3] \times \mathbb{C}/\mathbb{Z}[\zeta_3]\) with entropy the logarithm of a degree-four Salem number. The complete picture for entropies of projective torus automorphisms remains to be determined.

\section{K3 Surfaces}

A K3 surface is a simply connected two-dimensional compact complex manifold with a nowhere-vanishing holomorphic (2,0)-form. If \(Y\) is a K3 surface, then the intersection form makes \(H^2(Y;\mathbb{Z})\) into an even unimodular lattice of signature (3,19). As with two-dimensional complex tori, the induced action on second cohomology determines the entropy of any K3 surface automorphism. Indeed, the entropy of any K3 surface automorphism must be zero or the logarithm of a Salem number of degree at most twenty-two. (If the K3 surface is projective, then the degree of the Salem number cannot be greater than twenty.) (See \cite{mc1}, \S 3.)
\\
\\
A Kummer surface is a K3 surface constructed from a complex torus: let \(X\) be a two-dimensional complex torus and let \(i\) be the involution on \(X\) sending \(x\) to \(-x\); then the Kummer surface associated to \(X\) is the quotient \(X/i\) with the sixteen singularites resolved by blowing up; the blown-up quotient is simply-connected, and it inherits a nowhere-vanishing (2,0)-form from \(X\). If \(F\) is a torus automorphism, then \(F\) descends to an automorphism of the associated Kummer surface with the same entropy. Thus any Salem number that gives the entropy of some torus automorphism also gives the entropy of some K3 surface automorphism. A Kummer surface is projective if and only if it is constructed from a projective torus. If a Salem number gives the entropy of a Kummer surface automorphism that is induced (by the Kummer construction) from a torus automorphism, then the degree of the Salem number cannot be greater than six (and cannot be greater than four if the Kummer surface is projective). (See \cite{mc1}, \S 4.)
\\
\\
The following existence theorem leads to one method for identifying positive values that can be realized as entropies of K3 surface automorphisms. 
\begin{theorem}
[\cite{mc2}, Theorem 6.2] Let \(L_{3,19}\) be the unique even unimodular lattice of signature \((3,19)\). Let \(f\) be an isometry of \(L_{3,19}\) with spectral radius \(\rho(f)>1\) such that \(\rho(f)\) is an eigenvalue of \(f\). Suppose that there is some \(\tau \in (-2,2)\) such that
\[E_\tau=\ker(f+f^{-1}-\tau I) \subseteq L_{3,19} \otimes \mathbb{R}\]
has signature \((2,0)\). Then \(f\) is realizable by a K3 surface automorphism if \(f\) is the identity on \(L_{3,19} \cap E_\tau^\perp\).
\end{theorem}
Let \(L_{3,11}\) and \(L_{3,19}\) be the unique even unimodular lattices of signatures (3,11) and (3,19), respectively. Let \(S(t)\) be the minimal polynomial of a Salem number of degree fourteen or twenty-two. As in the degree-six case, results in \cite{gmc} show that \(S(t)\) is the characteristic polynomial of some isometry of \(L_{3,11}\) or \(L_{3,19}\) (depending on the degree) if \(S(-1)S(1)=-1\). In the degree-twenty-two case, this condition guarantees immediately that the Salem number gives the entropy of some K3 surface automorphism; this condition is also essential to the proof of Theorem 1.2. (Compare \cite{gmc}, \cite{mc2}, and \cite{ogu}.)
\\
\\
\emph{Proof of Theorem 1.2.} Let \(f\) be the isometry of \(L_{3,11}\) with characteristic polynomial \(S(t)\), let \(R(t)\) be the trace polynomial for \(S(t)\) given by \(S(t)=t^7R(t+t^{-1})\), and let \(E_8(-1)= \textless e_1, \ldots ,e_8 \textgreater\) be the even unimodular lattice of signature (0,8) whose intersection form is given by
\[(\langle e_i,e_j \rangle)=
\left( \begin{array}{cccccccc}
-2 & 1 & 0 & 0 & 0 & 0 & 0 & 0 \\
1 & -2 & 1 & 0 & 0 & 0 & 0 & 0 \\
0 & 1 & -2 & 1 & 1 & 0 & 0 & 0 \\
0 & 0 & 1 & -2 & 0 & 0 & 0 & 0 \\
0 & 0 & 1 & 0 & -2 & 1 & 0 & 0 \\
0 & 0 & 0 & 0 & 1 & -2 & 1 & 0 \\
0 & 0 & 0 & 0 & 0 & 1 & -2 & 1 \\
0 & 0 & 0 & 0 & 0 & 0 & 1 & -2 \end{array} \right).\] 
Then \(f\) extends to an isometry of \(L_{3,11} \oplus E_8(-1) \cong L_{3,19}\) which is the identity on the second summand, and there is a unique root \(\tau\) of \(R(t)\) giving \(E_\tau\) with signature (2,0). Since \(R(t)\) is irreducible over \(\mathbb{Z}\) and \(f+f^{-1}\) preserves \(E_\tau\), \(E_\tau^\perp\), and \(L_{3,11}\), it follows that
\[(L_{3,11} \oplus E_8(-1)) \cap E_\tau^\perp \subseteq E_8(-1),\]
on which \(f\) is the identity. \(\Box\)
\\
\\
Let \(L_{1,9}\) be the unique even unimodular lattice of signature (1,9) and let \(S(t)\) be the minimal polynomial of a Salem number of degree ten. The results in \cite{gmc} show again that \(S(t)\) is the characteristic polynomial of some isometry of \(L_{1,9}\) if \(S(-1)S(1)=-1\). However, in this case, no root of the trace polynomial \(R(t)\) given by \(S(t)=t^5R(t+t^{-1})\) gives rise to a space \(E_\tau\) with signature (2,0). Thus the approach used in the degree-twenty-two and degree-fourteen cases is not immediately applicable. In \cite{mc2}, there is a construction of a specific (non-projective) K3 surface automorphism with entropy the logarithm of a Salem number of degree ten (the smallest Salem number of degree ten, in fact). The construction relies on careful twisting and gluing of \(L_{1,9}\). Improvements to this process, in \cite{mc2} and \cite{mc3}, give constructions of projective K3 surface automorphisms whose entropies include logarithms of Salem numbers of degree six, eight, ten, and eighteen.

\section*{Acknowledgements}

We thank Laura DeMarco, Curt McMullen, and the referee for helpful comments.

\bibliographystyle{plain}
\bibliography{ReschkeP-refs-2012.02.22}

\end{document}